%% file: main.tex
\documentclass{amsart}
\makeatletter

\usepackage{import}
\import{./}{preamble}

\input lst-Macaulay2.tex
\definecolor{RED}{RGB}{255, 0, 0}

\title{A Step towards Computational Derived Algebraic Geometry : The RepHomology Package For Macaulay2}
\author{Guanyu Li}
\date{\today}

\begin{document}

\begin{abstract}
We introduce the \verb|Macaulay2| package \verb|RepHomology| for the computations of representation homology of certain spaces.
The main methods implement computing the representation homology of surfaces (with group coefficients, and analogies with algebra and Lie algebra coefficients), and the representation homology of link complements.
\end{abstract}

\maketitle

\section{Introduction}

Let $X$ be a ``nice'' connected topological space (e.g. CW complexes), and let $x_0$ be a fixed point in $X$. Then for any algebraic group $G$, the $G$-local systems (locally constant sheaves with a $G$-action) on $X$ is in a 1-1 correspondence of $\pi_1(X,x_0)$-representations in $G$.
It is known that all representations of $\pi_1(X,x_0)$ in $G$ are parameterized by an affine scheme $\mathrm{Rep}_G(\pi_1(X,x_0))$ which is called the representation scheme. In particular, the representation scheme also parameterizes all $G$-local systems on $X$.

These schemes are well studied, and some special cases, for instance the commuting scheme of a group scheme, connects many different areas of mathematics.
However, representation schemes do not generally behave well.
First, these schemes are generally very singular, and hence one needs to resolve the singularities of $\mathrm{Rep}_G(\pi_1(X,x_0))$.
Second, the scheme $\mathrm{Rep}_G(\pi_1(X,x_0))$ only relies on first order information $\pi_1(X,x_0)$ of the space $X$, and one needs to take into account higher homotopy information.

Derived algebraic geometry (DAG) offers a possible solution to remedy the deficiencies.
Kapranov \cite{K01} constructed the derived local system, which is a DG scheme relying on a certain injective resolution of $BG$ a priori.
Using the language of derived stacks by To\"en and Vezzosi \cite{TV08}, one has the derived mapping stack construction ${\bf {Map}}((X,x_0),(BG,*))$, whose underlying stack is the representation scheme $\mathrm{Rep}_G(\pi_1(X,x_0))$.
As the (nonabelian) derived functor of the representation scheme, the derived representation scheme is introduced by Berest, Ramadoss, and Yeung \cite{BRY22}.

All these derived objects are shown to be equivalent (\cite{BRY22-2}*{Appendix}), and the representation homology is the homology groups of the associated complexes to these derived objects.

The complexity of the object makes it difficult to calculate the homology groups.
There are very few known computations.
To facilitate further experimentation and exploration, we present a new \verb|Macaulay2| \cite{M2} package, \verb|RepHomology|.
This package aims to compute the representation homology of certain spaces, including surfaces, links, product of spheres, and etc.
There are also analogies to the representation homology of surfaces and product of spheres, via Lie algebra representation homology and associative algebra representation homology, and we give direct computations for these as well.

The code for this package can be found at
\href{https://github.com/GuanyuLee/RepresentationHomology}{https://github.com/GuanyuLee/RepresentationHomology}
and is expected to appear in a future release of \verb|Macaulay2|.

\subsection{Acknowledgements}

I want to express my greatest gratitude to Mike Stillman, who helped me a lot for this project.
I am also grateful to my advisor Yuri Berest, who introduced me this topic and discussed not only the big picture but also the details of computations.

\section{The complexes to compute representation homology}

Representation homology is the homology group of a certain complex associated to derived representation scheme.
We will explain the construction of the complexes here and give precise definitions in \Cref{Sec_RepH}.
In this paper, $k$ will always denote a field.

By a DG algebras, we mean a positively graded algebra with a differential of degree $-1$, satisfying the Leibnitz rule.
In particular, the differential of a $0$-degree element is always $0$.

\begin{eg}\label{Example_Running}
Let $X,Y$ be two generic upper triangular unipotent matrices over $k$
\begin{displaymath}
X=\begin{pmatrix}1&x_{1,2}&x_{1,3}&x_{1,4}\\
&1&x_{2,3}&x_{2,4}\\
&&1&x_{3,4}\\
&&&1
\end{pmatrix},\;\;\;\;
Y=\begin{pmatrix}1&y_{1,2}&y_{1,3}&y_{1,4}\\
&1&y_{2,3}&y_{2,4}\\
&&1&y_{3,4}\\
&&&1
\end{pmatrix},
\end{displaymath}
then $X,Y$ are generic matrices of the group scheme $U_4$, where $U_4$ consists of all $4\times4$ upper triangular unipotent matrices.
One sees that
\begin{equation}\label{Eq_CSU4Mat}
f:=XYX^{-1}Y^{-1}-I=\begin{pmatrix}
0&0&x_{1,2}y_{2,3}-x_{2,3}y_{1,2}&\substack{y_{1,2}x_{2,3}x_{3,4}-x_{1,2}y_{2,3}x_{3,4}\\ +y_{1,2}x_{2,3}y_{3,4}-x_{1,2}y_{2,3}y_{3,4}\\ +x_{1,2}y_{2,4}-y_{1,2}x_{2,4}+x_{1,3}y_{3,4}-y_{1,3}x_{3,4}}\\
&0&0&x_{2,3}y_{3,4}-x_{3,4}y_{2,3}\\
&&0&0\\
&&&0
\end{pmatrix}.
\end{equation}
Then the commuting scheme $C(U_4)$ of $U_4$ is
\begin{equation}
\Spec\frac{k[x_{i,j},y_{i,j}]}{f_{i,j}},\;\;1\leq i<j\leq4,
\end{equation}
where $f_{i,j}$ is the $(i,j)$-entry of the matrix $f$ defined in \Cref{Eq_CSU4Mat}.

Let $C$ be the Koszul (chain) complex of $\{f_{i,j}\}$ in $k[x_{i,j},y_{i,j}]$ ($1\leq i<j\leq4$).
Then the complex $C$ is an example computing the representation homology
\begin{displaymath}
HR_*(T^2,U_4)\cong H_*(C),
\end{displaymath}
where $T^2$ is the topological torus, and one can see that the $0$-th homology is exactly the coordinate ring of $C(U_4)$.
\end{eg}

In practice as in \Cref{Example_Running}, the algebraic groups (affine group schemes) are taken to be matrix groups, namely the groups consist of matrices and the multiplications come from matrix multiplication.
If $X$ is the generic matrix of the algebraic group $G$, then we denote by $k[X]$ the coordinate ring of $G$, i.e. $G=\Spec k[X]$.

In this package, we shall consider algebraic groups $GL_n,\,SL_n,\,U_n$ and $B_n$ (group consisting of upper triangular matrices with determinant invertible) in the package.
All these algebraic groups satisfy that the map $\epsilon:k[X]\to k$, which is induced by taking the identity in the group $G$, has a kernel generated by a regular sequence $\{z_1,\cdots,z_r\}$.

For an arbitrary topological space $X$, the complex compute the representation homology is generally infinite and difficult to manipulate.
However, for orientable compact surfaces and link complements in $\mathbb{R}^3$, there exist finite complexes to compute their representation homology, and they are stated below.

~\par
For surfaces, we have

\begin{prop}[See also {\cite{li2024commuting}*{Lemma 3.3}}]\label{Prop_SurfaceComplex}
Let $\Sigma_g$ be the compact orientable surface with genus $g$ and let $G$ be an affine matrix group scheme over $k$ with generic matrix $X$, whose coordinate ring is denoted by $k[X]$.
Then the DG algebra
\begin{equation}\label{Eq_SurfaceComplex}
k[X_1,\cdots,X_g,Y_1,\cdots,Y_g,T;\deg X_i=\deg Y_i=0,\deg T=1;dT=X_1Y_1X_1^{-1}Y_1^{-1}\cdots X_gY_gX_g^{-1}Y_g^{-1}-I]
\end{equation}
has homology groups isomorphic to $HR_*(\Sigma_g,G)$, where the matrix equation
\begin{equation*}
dT=X_1Y_1X_1^{-1}Y_1^{-1}\cdots X_gY_gX_g^{-1}Y_g^{-1}-I
\end{equation*}
means termwise equal for two matrices.
\end{prop}

And similarly for link complements in $\mathbb{R}^3$, we have

\begin{prop}[See also \cite{BRY22}*{Theorem 6.1}]\label{Prop_LinkComplex}
Let $L$ a link with braid word $\beta=\sigma_1\cdots\sigma_k\in B_n$ (i.e. $L=\hat{\beta}$, see \Cref{Sec_Link} for explanation), and let $G$ be an affine matrix group scheme over $k$ with generic matrix $X$, whose coordinate ring is denoted by $k[X]$.
Then the DG algebra
\begin{equation}
k[X_1,\cdots,X_n,Y_1,\cdots,Y_{n},T_1,\cdots,T_n;\deg X_i=\deg Y_i=0,\deg T_i=1;dT_i=Y_i-\beta_*(X_i)],
\end{equation}
has homology groups isomorphic to $HR_*(\Sigma_g,G)$, where $\beta_*$ is the composite operator $(\sigma_1)_*\cdots(\sigma_k)_*$ where $\sigma_i\in B_n$ and $(\sigma_i)_*$ is the operator mapping
\begin{displaymath}
X_j\mapsto\left\{\begin{matrix}
X_iX_{i+1}X_i^{-1},&\text{if }j=i\\
X_i,&\text{if }j=i+1\\
X_j,&\text{otherwise.}
\end{matrix}\right.
\end{displaymath}
\end{prop}

\begin{rmk}
One easily notices that in \Cref{Prop_SurfaceComplex} and \Cref{Prop_LinkComplex}, the DG algebras computing representation homology have underlying complexes being Koszul complexes.
This is simply a coincidence, and in general, one cannot expect the representation homology to be computed by a Koszul complex. In fact, within the same setting, for $S^m\times S^l$, where $m,l\geq3$ are odd natural numbers, the complex to compute $HR_*(S^m\times S^l,G)$ is
\begin{equation*}
k[X,Y,T;dT=[X,Y],\;\deg X=m,\;\deg Y=l,\;\deg T = m+l+1]
\end{equation*}
and it is not a Koszul complex simply by degrees.
Another such an example is the complex to compute $HR_*(L(p,q),G)$ where $L(p,q)$ is the lens space.
We hope to add the computation of $HR_*(S^m\times S^l,G)$ in a future release of our package.
\end{rmk}

\section{Representation Homology of Surfaces}

\subsection{Representations in groups}

Let $\Sigma_g$ be the compact orientable surface of genus $g$ as before, then one of the main goals of the package is the computation of
\begin{equation}\label{Eq_SurfaceGroup}
HR_*(\Sigma_g,G),
\end{equation}
which is the homology of the complex \Cref{Eq_SurfaceComplex}.
The interface function {\tt surfaceRepHomologyGroup} allows one to do the computation directly:

\begin{eg}\label{Example_FirstExample}
The following codes gives the representation homology $HR_*(\Sigma_1,U_3)$:
\begin{lstlisting}[language=Macaulay2]
`\underline{\tt i1}` : needsPackage "RepHomology";
`\underline{\tt i2}` : surfaceRepHomologyGroup(3, 1, GroupType=>"Unipotent")
H_0 = cokernel | x_(1,2,3)y_(1,1,2)-x_(1,1,2)y_(1,2,3) |
H_1 = cokernel | x_(1,2,3)y_(1,1,2)-x_(1,1,2)y_(1,2,3) 0                                     |
               | 0                                     x_(1,2,3)y_(1,1,2)-x_(1,1,2)y_(1,2,3) |
H_2 = cokernel | x_(1,2,3)y_(1,1,2)-x_(1,1,2)y_(1,2,3) |
`\underline{\tt o2}` = `...`
`\underline{\tt o2}` : List
\end{lstlisting}
One could compute the representation homology with other different groups. The default group type is $GL_n$:
\begin{lstlisting}[language=Macaulay2]
`\underline{\tt i3}` : surfaceRepHomologyGroup(2, 1)
H_0 = cokernel | x_(1,1,1)x_(1,2,2)y_(1,1,1)y_(1,2,1)s_1t_1+       ...                        |
H_1 = cokernel | x_(1,1,1)x_(1,2,2)y_(1,2,1)s_1t_1+        ...                                |
               | -x_(1,1,2)y_(1,2,1)^2s_1t_1-         ...                                     |
               | x_(1,1,1)y_(1,2,1)s_1t_1+            ...                                     |
H_2 = cokernel | 0                       -x_(1,2,1)s_1t_1   ...                               |
               | x_(1,1,1)x_(1,2,1)s_1^2 x_(1,1,1)x_(1,2,2)y_(1,2,1)s_1^2t_1+...              |
               | -1                      y_(1,1,1)t_1-2y_(1,2,2)t_1   ...                     |
               | -x_(1,2,1)s_1         -x_(1,2,2)y_(1,2,1)s_1t_1-x_(1,2,1)y_(1,2,2)s_1t_1 ... |
`\underline{\tt o3}` = `...`
`\underline{\tt o3}` : List
\end{lstlisting}
We can also take the \verb|GroupType| to be \verb|"SL"| or \verb|"Borel"|.
\end{eg}

In particular, we could also verify \cite{li2024commuting}*{Corollary 3.2.1}:

\begin{eg}
We would like to check
when $n=2,3,4,5$ the condition
\begin{equation}\label{Eq_Vanishing}
HR_i(\Sigma_1,U_n)=0\;\;\;\forall i\geq n.
\end{equation}

\begin{lstlisting}[language=Macaulay2]
`\underline{\tt i4}` : surfaceRepHomologyGroup(2, 1, GroupType=>"Unipotent")
                          1
H_0 = (QQ[x     , y     ])
           1,1,2   1,1,2
                          1
H_1 = (QQ[x     , y     ])
           1,1,2   1,1,2
`\underline{\tt o4}` = `...`
`\underline{\tt o4}` : List
`\underline{\tt i5}` : surfaceRepHomologyGroup(4, 1, GroupType=>"Unipotent")
H_0 = cokernel | x_(1,3,4)y_(1,2,3)-x_(1,2,3)y_(1,3,4) ... |
H_1 = cokernel | x_(1,3,4)y_(1,2,3)-x_(1,2,3)y_(1,3,4) ... |
               | 0                                     ... |
               | 0                                     ... |
H_2 = cokernel | x_(1,3,4)y_(1,2,3)-x_(1,2,3)y_(1,3,4) ... |
               | 0                                     ... |
               | 0                                     ... |
H_3 = cokernel | x_(1,3,4)y_(1,2,3)-x_(1,2,3)y_(1,3,4) ... |
`\underline{\tt o5}` = `...`
`\underline{\tt o5}` : List
`\underline{\tt i6}` : surfaceRepHomologyGroup(5, 1, GroupType=>"Unipotent")
H_0 = cokernel | x_(1,4,5)y_(1,3,4)-x_(1,3,4)y_(1,4,5) ... |
H_1 = cokernel | x_(1,4,5)y_(1,3,4)-x_(1,3,4)y_(1,4,5) ... |
               | 0                                     ... |
               | 0                                     ... |
               | 0                                     ... |
H_2 = cokernel | x_(1,4,5)y_(1,3,4)-x_(1,3,4)y_(1,4,5) ... |
               | 0                                     ... |
               | 0                                     ... |
               | 0                                     ... |
               | 0                                     ... |
               | 0                                     ... |
H_3 = cokernel | x_(1,4,5)y_(1,3,4)-x_(1,3,4)y_(1,4,5) ... |
               | 0                                     ... |
               | 0                                     ... |
               | 0                                     ... |
H_4 = cokernel | x_(1,4,5)y_(1,3,4)-x_(1,3,4)y_(1,4,5) ... |
`\underline{\tt o6}` = `...`
`\underline{\tt o6}` : List
\end{lstlisting}
Together with \Cref{Example_FirstExample}, we verify the condition (\ref{Eq_Vanishing}).
\end{eg}

\subsection{Representations in algebras}

Given a finite dimensional associative $k$-algebra $B$, for instance
\begin{displaymath}
B=N_m=\{A\in\mathrm{Mat}_{m\times m}(k)\mid A\text{ is strict upper triangular}\}.
\end{displaymath}
with usual matrix multiplication.
Take the associative algebra
\begin{displaymath}
\Pi_0(Q_g):=\frac{k<x_1,y_1,\cdots,x_g,y_g>}{<x_1y_1-y_1x_1+\cdots+x_gy_g-y_gx_g>}
\end{displaymath}
where $k<x_1,y_1,\cdots,x_n,y_n>$ is the free associative algebra generated by $x_1,y_1,\cdots,x_g,y_g$ and $<x_1y_1-y_1x_1+\cdots+x_gy_g-y_gx_g>$ is the 2-sided ideal generated by a single element
\begin{displaymath}
x_1y_1-y_1x_1+\cdots+x_gy_g-y_gx_g.
\end{displaymath}
Then the representation homology
\begin{equation}\label{Eq_SurfaceAlg}
HR_*(\Pi_0(Q_g),B)
\end{equation}
is an analog of \Cref{Eq_SurfaceGroup}.
The DG algebra to compute \Cref{Eq_SurfaceAlg} is
\begin{equation}\label{Eq_SurfaceAlgComplex}
k[X_1,\cdots,X_g,Y_1,\cdots,Y_g,T;\deg X_i=\deg Y_i=0,\deg T=1;dT=[X_1,Y_1]+\cdots+[X_g,Y_g]]
\end{equation}
where $[X_i,Y_i]=X_iY_i-Y_iX_i$.
We have the function {\tt surfaceRepHomologyAlg} to do the computation directly:

\begin{eg}\label{Example_AlgRepH}
The following codes gives the representation homology $HR_*(\Pi_0(Q_1),N_3)$, where $N_3$ is the matrix algebra consisting of all $3\times3$ strict upper triangular matrices.
\begin{lstlisting}[language=Macaulay2]
`\underline{\tt i7}` : surfaceRepHomologyAlg(3, 1, AlgType=>"nilpotent")
H_0 = cokernel | x_(1,2,3)y_(1,1,2)-x_(1,1,2)y_(1,2,3) |
H_1 = cokernel | x_(1,2,3)y_(1,1,2)-x_(1,1,2)y_(1,2,3) 0                                     |
               | 0                                     x_(1,2,3)y_(1,1,2)-x_(1,1,2)y_(1,2,3) |
H_2 = cokernel | x_(1,2,3)y_(1,1,2)-x_(1,1,2)y_(1,2,3) |
`\underline{\tt o7}` = `...`
`\underline{\tt o7}` : List
\end{lstlisting}

From the computational results in \Cref{Example_FirstExample}, one sees that there is an isomorphism
\begin{displaymath}
HR_*(\Sigma_1,U_3)\cong HR_*(\Pi_0(Q_1),N_3),
\end{displaymath}
and this is not a coincidence.
\end{eg}

Similarly to \Cref{Example_AlgRepH}, we can also take the \verb|AlgType| to be \verb|"gl"|, \verb|"sl"|, \verb|"nilpotent"|, or \verb|"borel"|.

\subsection{Representations in Lie algebras}

Let $\mathfrak{g}$ be a finite dimensional Lie algebra, for instance $\mathfrak{gl}_n$. Let $\mathfrak{a}$ be the $2$-dimensional abelian Lie algebra.
There is a function \verb|surfaceRepHomologyLie| to compute
\begin{equation}\label{Eq_TorusLie}
HR_*(\mathfrak{a},\mathfrak{g})
\end{equation}
directly. The complex to compute \Cref{Eq_TorusLie} is
\begin{equation}\label{Eq_SurfaceLieComplex}
k[X,Y,T;\deg X=\deg Y=0,\deg T=1;dT=[X,Y]]
\end{equation}
This is a special case of \Cref{Eq_SurfaceAlgComplex}. And this phenomenon can be generalised to broader cases.

\section{Representation Homology of Links}

\subsection{The type {\tt Link}}\label{Sec_Link}

This package also includes a new type \verb|Link| and basic operations of \verb|Link| objects.

We collect basic facts about links here first.

By definition, a (topological) link $L$ is a smooth oriented embedding of the disjoint union $S^1\coprod S^1\coprod\cdots\coprod S^1$ into $\mathbb{R}^3$.
The link complement $X:=\mathbb{R}^3-L$ is defined to be the complement of an open tubular neighborhood of the image of $L$ in $\mathbb{R}^3$.

It is known that every link $L$ in $\mathbb{R}^3$ can be obtained geometrically as the closure of a braid $\beta\in B_n$ in $\mathbb{R}^3$ (we write $L=\hat{\beta}$ to indicate this relation) \cite{B74}.
Thus one could use a braid to represent a link complement in $\mathbb{R}^3$.

Link is a new type of \verb|HashTable|.
A link consists of a \verb|BraidIndex| and a \verb|LinkWord|, where the \verb|BraidIndex| is an integer $n$ indicating the size of braid group $B_n$, and \verb|LinkWord| is a list of nonzero integers indicating the braid operations.
For instance, if
\begin{displaymath}
\beta=\sigma_1^3\in B_2
\end{displaymath}
then $\hat{\beta}$ is the simplest nontrivial knot and one could construct this link as follows
\begin{lstlisting}[language=Macaulay2]
`\underline{\tt i8}` : L = link(2,{1,1,1})
`\underline{\tt o8}` = L
`\underline{\tt o8}` : Link
\end{lstlisting}
where the integer $2$ means the braid is taken from $B_2$ and the list $\{1,1,1\}$ means the braid operation is $\sigma_1\circ\sigma_1\circ\sigma_1=\sigma_1^3$. For the inverse operation $\sigma_i^{-1}$, we use $-i$ in the list to denote it.

There are two built-in links
\begin{lstlisting}[language=Macaulay2]
`\underline{\tt i9}` : Trefoil
`\underline{\tt o9}` = Trefoil
`\underline{\tt o9}` : Link
`\underline{\tt i10}` : FigureEight
`\underline{\tt o10}` = FigureEight
`\underline{\tt o10}` : Link
\end{lstlisting}
which are $\sigma_1^3$ and $\sigma_1\sigma_2^{-1}\sigma_1\sigma_2^{-1}$ using braid words.

One could also see \cite{J87} for a list of knots by braid words.

\subsection{Representation Homology of Link Complements}

Given a link $L$, one could compute the representation homology of the link complement in $\mathbb{R}^3$ directly by the function {\tt linkRepHomology}.

\begin{eg}
The following codes gives the representation homology $HR_*(L,U_3)$ where $L$ is the knot trefoil:
\begin{lstlisting}[language=Macaulay2]
`\underline{\tt i11}` : linkRepHomology(Trefoil,3,GroupType=>"U")
H_0 = cokernel | y_(1,2,3)-y_(2,2,3) y_(1,1,3)-y_(2,1,3) ... |
H_1 = cokernel | ... |
               | ... |
               | ... |
H_2 = cokernel | ... |
               | ... |
               | ... |
H_3 = cokernel | -y_(1,2,3)+y_(2,2,3) -y_(1,1,3)+y_(2,1,3) ... |
`\underline{\tt o11}` = `...`
`\underline{\tt o11}` : List
\end{lstlisting}

\end{eg}

\appendix

\section{Representation homology}\label{Sec_RepH}

We will give a short introduction to the precise definition of representation homology (using Berest-Ramadoss-Yeung's derived representation scheme model) here, as well as some analogies.

\subsection{Notations and conventions}

We denote by $\Top$ and $\Top_*$ the categories of nice topological spaces and nice pointed topological spaces, by $\Gp$ the category of groups, and $\Sch_k$ the category of schemes over $\Spec k$.

The category of all DG commutative (resp. associative and Lie) algebras will be denoted by $\CDGA_k$ (resp. ${\bf DGA}_k$ and ${\bf DGLA}_k$).

\subsection{Simplicial methods}

We denote by $\bf\Delta$ the category of simplices (see \cite{GJ09}).

Given any category $\cC$, denote by $s\cC$ the category of simplicial objects of $\cC$, i.e. $s\cC:=\mathrm{Funct}(\bf\Delta^{\mathrm{op}},\cC)$.
A simplicial set $X$ is called reduced if $X_0=\{*\}$. The full subcategory of $s\Set$ consisting of reduced simplicial sets will be
denoted $s\Set_0$.

There is a Quillen equivalence
\begin{displaymath}
|-|:s\Set\rightleftarrows\Top:S
\end{displaymath}
where $S$ is the (total) singular complex and $|-|$ is the geometric realisation (\cite{GJ09}*{Theorem I.11.4.}).

Restricted to the category $\Top_{0,*}$ of connected pointed spaces, we get the pair of adjoint functors
\begin{displaymath}
|-|:s\Set_0\rightleftarrows\Top_{0,*}:\bar{S}
\end{displaymath}
where $\bar{S}_*(X)$ is defined by
\begin{displaymath}
\bar{S}_n(X):=\{f:\Delta^n\to X\mid f(v_i)=* \text{ for all vertices }v_i\in\Delta^n\}.
\end{displaymath}
~\par

Furthermore, there is another Quillen equivalence
\begin{equation}\label{Eq_KanLoopGroup}
\mathbb{G}:s\Set_0\rightleftarrows s\Gp:\bar{W},
\end{equation}
where $\mathbb{G}$ is called the Kan loop group functor and $\bar{W}$ is the classifying simplicial complex (see {\cite{GJ09}*{Proposition V.6.3.}}). It is worthwhile to point out that by the construction, the Kan loop group $\mathbb{G}X$ of any reduced simplicial set $X$ is semi-free, hence cofibrant. In this way, one could have a ``resolution'' $\mathbb{G}(X)$ of the space $X$.

In conclusion we have a sequence of equivalences of homotopy categories
\begin{equation}\label{Eq_Simplicial}
\Ho(\Top_{0,*})\simeq\Ho(s\Set_0)\simeq\Ho(s\Gp).
\end{equation}
This means giving a simplicial group is the same as giving a pointed, connected topological group.

The famous Dold-Kan correspondence states that there is a Quillen equivalence
\begin{displaymath}
N:s\Vect_k\rightleftarrows\Com(\Vect_k):\Gamma
\end{displaymath}
when $N$ is the normalisation functor and $\Gamma$ is its right adjoint (\cite{GJ09}*{Corollary III.2.3, Corollary III.2.7.}).
Then for any simplicial vector space $V\in s\Vect_k$, we define the homotopy group $\pi_*(V)$ to be the homology group of its normalisation $N(V)$, i.e.
\begin{equation*}
\pi_*(V):=H_*(N(V)).
\end{equation*}
~\par

\subsection{Representation homology}

Given an affine group scheme $G$ over $k$, by viewing it a functor (of points), one has an adjunction
\begin{equation}\label{Eq_BasicAdj}
(-)_G:\Gp\leftrightarrows{\bf CommAlg}_k:G
\end{equation}
where the left functor $(-)_G:\Gp\to{\bf CommAlg}_k$ exists since a group homomorphism only satisfies algebraic equations. To any group $\Gamma$, the scheme $\mathrm{Rep}_G(\Gamma):=\Spec(\Gamma)_G$ associated to the ring $(\Gamma)_G$ is called the \textbf{representation scheme}.
This adjunction \ref{Eq_BasicAdj} could be extended simplicially, getting another adjunction
\begin{equation}\label{Eq_ExtendedAdj}
(-)_G:s\Gp\leftrightarrows s{\bf CommAlg}_k:G.
\end{equation}
This is not a Quillen pair, but we still have

\begin{lem}[{\cite{BRY22}*{Lemma 3.1}}]\label{Lemma_ExistenceOfDerivedFunctor}
The functor $(-)_G$ in (\ref{Eq_ExtendedAdj}) maps the weak equivalences between cofibrant objects in $s\Gp$ to weak equivalences in $s{\bf CommAlg}_k$, and hence has a total left derived functor
\begin{equation}\label{Eq_DerivedRepSch}
\mathbb{L}(-)_G:\Ho(s\Gp)\to\Ho(s\Alg_k).
\end{equation}
As a consequence, the left derived functor $\mathbb{L}(-)_G$ can be computed by $X\mapsto(QX)_G$, where $QX$ is a cofibrant replacement of $X$ in $s\Gp$.
\end{lem}

We shall denote the derived functor by $\mathscr{O}(\mathrm{DRep}_G(\Gamma)):\mathrm{Ho}(s\Gp)\to\mathrm{Ho}(s\Alg_k)$ and call the corresponding (derived) scheme $\mathrm{DRep}_G(\Gamma)$ the {\bf derived representation scheme}.

\begin{defn}
Given a (pointed, connected) topological space $X$ and a group scheme $X$ over $k$, the representation homology $HR_*(X,G)$ of $X$ with coefficient $G$ is the homotopy group
\begin{displaymath}
\pi_*((\mathbb{G}(X))_G).
\end{displaymath}
\end{defn}

One should notice that one of the key implications of \Cref{Lemma_ExistenceOfDerivedFunctor} is that the representation homology is independent of the choice of the cofibrant replacement of $X$ in $s\Gp$.
Therefore we can take a certain nice approachable cofibrant replacement of the space $X$, instead of the standard Kan loop group construction $\mathbb{G}(X)$ which is huge and almost impossible to do computations.
The nice resolution does not always exist, but in many cases such as orientable surfaces and link complements, we could find a Koszul complex as a resolution, which allows us to calculate the homology.

The following theorem is the key to the calculations.

\begin{thm}[{\cite{BRY22}*{Theorem 3.2}}]
The derived representation functor \Cref{Eq_DerivedRepSch} preserves homotopy colimits.
\end{thm}

\subsection{Other types of representation homology}\label{Sec_Others}

We assume the characteristic of the field $k$ to be $0$ in this subsection.

There are analogies of the simplicially extended adjunction \Cref{Eq_ExtendedAdj}, using the representation schemes for associative algebras and the representation schemes for Lie algebras.
Precisely, given any finite dimensional (associative) algebra $B$ there is an adjunction
\begin{displaymath}
(-)_{B}:{\bf Alg}_k\leftrightarrows{\bf CommAlg}_k:\mathrm{Hom}_k(B^*,-)
\end{displaymath}
where $\mathrm{Hom}_k(B^*,-)$ gives the convolution algebra and $(-)_{B}$ is the representation scheme functor for associative algebras. This adjunction extends to
\begin{equation}\label{Eq_DRepAlg}
(-)_B:{\bf DGA}_k\to{\bf CDGA}_k:\mathrm{Hom}_k(B^*,-)
\end{equation}
which is a Quillen pair.
Thus the derived functor $\mathbb{L}(-)_{B}$ of $(-)_{B}$ exists and for an DG associative algebra $A$, we call
\begin{equation}
HR(A,B):=H_*(\mathbb{L}(A)_{B})
\end{equation}
the representation homology of $A$ in $B$.

Similarly, given any finite dimensional Lie algebra $\mathfrak{g}$ there is an adjunction
\begin{displaymath}
(-)_{\mathfrak{g}}:{\bf Lie}_k\to{\bf CommAlg}_k:\mathfrak{g}\otimes-
\end{displaymath}
where $(-)_{\mathfrak{g}}$ is the representation scheme functor for Lie algebras. This adjunction extends to
\begin{equation}\label{Eq_DRepLie}
(-)_{\mathfrak{g}}:{\bf DGLA}_k\to{\bf CDGA}_k:\mathfrak{g}\otimes-
\end{equation}
which is also a Quillen pair.
Thus the derived functor $\mathbb{L}(-)_{\mathfrak{g}}$ of $(-)_{\mathfrak{g}}$ exists and for an DG Lie algebra $\mathfrak{a}$, we call
\begin{equation}
HR(\mathfrak{a},\mathfrak{g}):=H_*(\mathbb{L}(\mathfrak{a})_{\mathfrak{g}})
\end{equation}
the representation homology of $\mathfrak{a}$ in $\mathfrak{g}$.

\end{document}

%% file: preamble.tex
\usepackage[T1]{fontenc}
\usepackage{textcmds}  
\usepackage{amsmath, amssymb, amsfonts, amstext, verbatim, amsthm, mathrsfs}
\numberwithin{equation}{section} 
\usepackage{geometry}
\usepackage{mathtools}
\usepackage[mathcal]{eucal}
\usepackage{microtype}
\usepackage{thmtools, thm-restate}
\usepackage[all]{xy}
\usepackage{tikz-cd}
\usepackage{tikz} 
\usepackage{enumerate}
\usepackage{enumitem}
\usepackage[modulo]{lineno}
\usepackage[colorlinks=true,linkcolor=blue,citecolor=blue,urlcolor=blue,citebordercolor={0 0 1},urlbordercolor={0 0 1},linkbordercolor={0 0 1}]{hyperref} 
\usepackage[shortalphabetic]{amsrefs} 
\usepackage[nameinlink]{cleveref}

\geometry{top=2.5cm,bottom=2.5cm,left=2.1cm,right=2.1cm}

\def\makeCal#1{%
\expandafter\newcommand\csname c#1\endcsname{\mathcal{#1}}}
\def\makeBB#1{%
\expandafter\newcommand\csname b#1\endcsname{\mathbb{#1}}}
\def\makeFrak#1{%
\expandafter\newcommand\csname f#1\endcsname{\mathfrak{#1}}}
\def\makeScr#1{%
\expandafter\newcommand\csname s#1\endcsname{\mathscr{#1}}}

\count@=0
\loop
\advance\count@ 1
\edef\y{\@Alph\count@}%
\expandafter\makeCal\y
\expandafter\makeBB\y
\expandafter\makeFrak\y
\ifnum\count@<26
\repeat


\newcommand{\Alg}{\mathbf{Alg}}

\newcommand{\CDGA}{\mathbf{CDGA}}

\newcommand{\Com}{\mathbf{Com}}

\newcommand{\Gp}{\mathbf{Gp}}

\DeclareMathOperator{\Ho}{Ho}

\newcommand{\Sch}{\mathbf{Sch}}
\newcommand{\Set}{\mathbf{Set}}

\DeclareMathOperator{\Spec}{Spec~}

\newcommand{\Top}{\mathbf{Top}}

\newcommand{\Vect}{\mathbf{Vect}}



\theoremstyle{plain}
\newtheorem{thm}{Theorem}[section]

\newtheorem{lem}[thm]{Lemma}

\newtheorem{prop}[thm]{Proposition}

\theoremstyle{definition}
\newtheorem{defn}{Definition}[section]

\newtheorem{eg}[thm]{Example}

\theoremstyle{remark}
\newtheorem{rmk}{Remark}

\AtEndDocument{%
\bibliographystyle{unsrt}
\nocite{*}
\bibliography{./reference.bib}
}

%% file: lst-Macaulay2.tex
\usepackage{listings}
\usepackage{xcolor}
\definecolor{maccolor}{rgb}{0.3,0.3,0.8}
\lstdefinelanguage{Macaulay2}
{
basicstyle={\ttfamily},
keywordstyle={\color{maccolor!80!black}},
commentstyle={\color{gray}},
stringstyle={\color{red!40!black}},
rulecolor=\color{maccolor},
basewidth={1.2ex}, 
sensitive=false,
morecomment=[l]{--},
morecomment=[s]{-*}{*-},
morestring=[b]",
escapechar={`},
escapebegin={\rmfamily},
morekeywords={about,abs,AbstractToricVarieties,accumulate,Acknowledgement,acos,acosh,acot,addCancelTask,addDependencyTask,addEndFunction,addHook,AdditionalPaths,addStartFunction,addStartTask,Adjacent,adjoint,AdjointIdeal,AffineVariety,AfterEval,AfterNoPrint,AfterPrint,agm,AInfinity,alarm,AlgebraicSplines,Algorithm,Alignment,all,AllCodimensions,allowableThreads,ambient,analyticSpread,Analyzer,AnalyzeSheafOnP1,ancestor,ancestors,ANCHOR,and,andP,AngleBarList,ann,annihilator,antipode,any,append,applicationDirectory,applicationDirectorySuffix,apply,applyKeys,applyPairs,applyTable,applyValues,apropos,argument,Array,arXiv,Ascending,ascii,asin,asinh,ass,assert,associatedGradedRing,associatedPrimes,AssociativeAlgebras,AssociativeExpression,atan,atan2,atEndOfFile,Authors,autoload,AuxiliaryFiles,backtrace,Bag,Bareiss,baseFilename,BaseFunction,baseName,baseRing,baseRings,BaseRow,BasicList,basis,BasisElementLimit,Bayer,BeforePrint,beginDocumentation,BeginningMacaulay2,Benchmark,benchmark,Bertini,BesselJ,BesselY,betti,BettiCharacters,BettiTally,between,BGG,BIBasis,Binary,BinaryOperation,Binomial,binomial,BinomialEdgeIdeals,Binomials,BKZ,BlockMatrix,BLOCKQUOTE,BODY,Body,BoijSoederberg,BOLD,Book3264Examples,Boolean,BooleanGB,borel,Boxes,BR,break,Browse,Bruns,cache,CacheExampleOutput,CacheFunction,CacheTable,cacheValue,CallLimit,cancelTask,capture,catch,Caveat,CC,CDATA,ceiling,Center,centerString,Certification,ChainComplex,chainComplex,ChainComplexExtras,ChainComplexMap,ChainComplexOperations,ChangeMatrix,char,CharacteristicClasses,characters,charAnalyzer,check,CheckDocumentation,chi,Chordal,class,Classic,clean,clearAll,clearEcho,clearOutput,close,closeIn,closeOut,ClosestFit,CODE,code,codim,CodimensionLimit,coefficient,CoefficientRing,coefficientRing,coefficients,Cofactor,CohenEngine,CohenTopLevel,CoherentSheaf,CohomCalg,cohomology,coimage,CoincidentRootLoci,coker,cokernel,collectGarbage,columnAdd,columnate,columnMult,columnPermute,columnRankProfile,columnSwap,combine,Command,commandInterpreter,commandLine,COMMENT,commonest,commonRing,comodule,CompactMatrix,compactMatrixForm,CompiledFunction,CompiledFunctionBody,CompiledFunctionClosure,Complement,complement,complete,CompleteIntersection,CompleteIntersectionResolutions,Complexes,ComplexField,components,compose,compositions,compress,concatenate,conductor,ConductorElement,cone,Configuration,ConformalBlocks,conjugate,connectionCount,Consequences,Constant,Constants,constParser,content,continue,contract,Contributors,ConvexInterface,conwayPolynomial,ConwayPolynomials,copy,copyDirectory,copyFile,copyright,Core,CorrespondenceScrolls,cos,cosh,cot,CotangentSchubert,cotangentSheaf,coth,cover,coverMap,cpuTime,createTask,Cremona,csc,csch,current,currentColumnNumber,currentDirectory,currentFileDirectory,currentFileName,currentLayout,currentLineNumber,currentPackage,currentString,currentTime,Cyclotomic,Database,Date,DD,dd,deadParser,debug,debugError,DebuggingMode,debuggingMode,debugLevel,DecomposableSparseSystems,Decompose,decompose,deepSplice,Default,default,defaultPrecision,Degree,degree,degreeLength,DegreeLift,DegreeLimit,DegreeMap,DegreeOrder,DegreeRank,Degrees,degrees,degreesMonoid,degreesRing,delete,demark,denominator,Dense,Density,Depth,depth,Descending,Descent,Describe,describe,Description,det,determinant,DeterminantalRepresentations,DGAlgebras,diagonalMatrix,diameter,Dictionary,dictionary,dictionaryPath,diff,DiffAlg,difference,dim,directSum,disassemble,discriminant,dismiss,Dispatch,distinguished,DIV,Divide,divideByVariable,DivideConquer,DividedPowers,Divisor,DL,Dmodules,do,doc,docExample,docTemplate,document,DocumentTag,Down,drop,DT,dual,eagonNorthcott,EagonResolution,echoOff,echoOn,EdgeIdeals,edit,EigenSolver,eigenvalues,eigenvectors,eint,EisenbudHunekeVasconcelos,elapsedTime,elapsedTiming,elements,Eliminate,eliminate,Elimination,EliminationMatrices,EllipticCurves,EllipticIntegrals,else,EM,Email,End,end,endl,endPackage,Engine,engineDebugLevel,EngineRing,EngineTests,entries,EnumerationCurves,environment,Equation,EquivariantGB,erase,erf,erfc,error,errorDepth,euler,EulerConstant,eulers,even,EXAMPLE,ExampleFiles,ExampleItem,examples,ExampleSystems,Exclude,exec,exit,exp,expectedReesIdeal,expm1,exponents,export,exportFrom,exportMutable,Expression,expression,Ext,extend,ExteriorIdeals,ExteriorModules,exteriorPower,Factor,factor,false,Fano,FastMinors,FastNonminimal,FGLM,File,fileDictionaries,fileExecutable,fileExists,fileExitHooks,fileLength,fileMode,FileName,FilePosition,fileReadable,fileTime,fileWritable,fillMatrix,findFiles,findHeft,FindOne,findProgram,findSynonyms,FiniteFittingIdeals,First,first,firstkey,FirstPackage,fittingIdeal,flagLookup,FlatMonoid,flatten,flattenRing,Flexible,flip,floor,flush,fold,FollowLinks,for,forceGB,fork,FormalGroupLaws,Format,format,formation,FourierMotzkin,FourTiTwo,fpLLL,frac,fraction,FractionField,frames,FrobeniusThresholds,from,fromDividedPowers,fromDual,Function,FunctionApplication,FunctionBody,functionBody,FunctionClosure,FunctionFieldDesingularization,fusePairs,futureParser,GaloisField,Gamma,gb,GBDegrees,gbRemove,gbSnapshot,gbTrace,gcd,gcdCoefficients,gcdLLL,GCstats,genera,GeneralOrderedMonoid,GenerateAssertions,generateAssertions,generator,generators,Generic,GenericInitialIdeal,genericMatrix,genericSkewMatrix,genericSymmetricMatrix,gens,genus,get,getc,getChangeMatrix,getenv,getGlobalSymbol,getNetFile,getNonUnit,getPrimeWithRootOfUnity,getSymbol,getWWW,GF,gfanInterface,Givens,GKMVarieties,GLex,Global,global,globalAssign,globalAssignFunction,GlobalAssignHook,globalAssignment,globalAssignmentHooks,GlobalDictionary,GlobalHookStore,globalReleaseFunction,GlobalReleaseHook,Gorenstein,GradedLieAlgebras,GradedModule,gradedModule,GradedModuleMap,gradedModuleMap,gramm,GraphicalModels,GraphicalModelsMLE,Graphics,graphIdeal,graphRing,Graphs,Grassmannian,GRevLex,GroebnerBasis,groebnerBasis,GroebnerBasisOptions,GroebnerStrata,GroebnerWalk,groupID,GroupLex,GroupRevLex,GTZ,Hadamard,handleInterrupts,HardDegreeLimit,hash,HashTable,hashTable,HEAD,HEADER1,HEADER2,HEADER3,HEADER4,HEADER5,HEADER6,HeaderType,Heading,Headline,Heft,heft,Height,height,help,Hermite,hermite,Hermitian,HH,hh,HigherCIOperators,HighestWeights,Hilbert,hilbertFunction,hilbertPolynomial,hilbertSeries,HodgeIntegrals,hold,Holder,Hom,homeDirectory,HomePage,Homogeneous,Homogeneous2,homogenize,homology,homomorphism,HomotopyLieAlgebra,hooks,horizontalJoin,HorizontalSpace,HR,HREF,HTML,html,httpHeaders,Hybrid,HyperplaneArrangements,Hypertext,hypertext,HypertextContainer,HypertextParagraph,icFracP,icFractions,icMap,icPIdeal,id,Ideal,ideal,idealizer,identity,if,IgnoreExampleErrors,ii,image,imaginaryPart,IMG,ImmutableType,importFrom,in,incomparable,Increment,independentSets,indeterminate,IndeterminateNumber,Index,index,indexComponents,IndexedVariable,IndexedVariableTable,indices,inducedMap,inducesWellDefinedMap,InexactField,InexactFieldFamily,InexactNumber,InfiniteNumber,infinity,info,InfoDirSection,infoHelp,Inhomogeneous,input,Inputs,insert,installAssignmentMethod,installedPackages,installHilbertFunction,installMethod,installMinprimes,installPackage,InstallPrefix,instance,instances,IntegralClosure,integralClosure,integrate,IntermediateMarkUpType,interpreterDepth,intersect,intersectInP,Intersection,intersection,interval,InvariantRing,inverse,InverseMethod,inversePermutation,Inverses,inverseSystem,InverseSystems,Invertible,InvolutiveBases,irreducibleCharacteristicSeries,irreducibleDecomposition,isAffineRing,isANumber,isBorel,isCanceled,isCommutative,isConstant,isDirectory,isDirectSum,isEmpty,isField,isFinite,isFinitePrimeField,isFreeModule,isGlobalSymbol,isHomogeneous,isIdeal,isInfinite,isInjective,isInputFile,isIsomorphism,isLinearType,isListener,isLLL,isMember,isModule,isMonomialIdeal,isNormal,isOpen,isOutputFile,isPolynomialRing,isPrimary,isPrime,isPrimitive,isPseudoprime,isQuotientModule,isQuotientOf,isQuotientRing,isReady,isReal,isReduction,isRegularFile,isRing,isSkewCommutative,isSorted,isSquareFree,isStandardGradedPolynomialRing,isSubmodule,isSubquotient,isSubset,isSupportedInZeroLocus,isSurjective,isTable,isUnit,isWellDefined,isWeylAlgebra,ITALIC,Iterate,Jacobian,jacobian,jacobianDual,Jets,Join,join,Jupyter,K3Carpets,K3Surfaces,Keep,KeepFiles,KeepZeroes,ker,kernel,kernelLLL,kernelOfLocalization,Key,keys,Keyword,Keywords,kill,koszul,Kronecker,KustinMiller,LABEL,last,lastMatch,LATER,LatticePolytopes,Layout,lcm,leadCoefficient,leadComponent,leadMonomial,leadTerm,Left,left,length,LengthLimit,letterParser,Lex,LexIdeals,LI,Licenses,LieTypes,lift,liftable,Limit,limitFiles,limitProcesses,Linear,LinearAlgebra,LinearTruncations,lineNumber,lines,LINK,linkFile,List,list,listForm,listLocalSymbols,listSymbols,listUserSymbols,LITERAL,LLL,LLLBases,lngamma,load,loadDepth,LoadDocumentation,loadedFiles,loadedPackages,loadPackage,Local,local,localDictionaries,LocalDictionary,localize,LocalRings,locate,log,log1p,LongPolynomial,lookup,lookupCount,LowerBound,LUdecomposition,M0nbar,M2CODE,Macaulay2Doc,makeDirectory,MakeDocumentation,makeDocumentTag,MakeHTML,MakeInfo,MakeLinks,makePackageIndex,MakePDF,makeS2,Manipulator,map,MapExpression,MapleInterface,markedGB,Markov,MarkUpType,match,mathML,Matrix,matrix,MatrixExpression,Matroids,max,maxAllowableThreads,maxExponent,MaximalRank,maxPosition,MaxReductionCount,MCMApproximations,member,memoize,memoizeClear,memoizeValues,MENU,merge,mergePairs,META,method,MethodFunction,MethodFunctionBinary,MethodFunctionSingle,MethodFunctionWithOptions,methodOptions,methods,midpoint,min,minExponent,mingens,mingle,minimalBetti,MinimalGenerators,MinimalMatrix,minimalPresentation,minimalPresentationMap,minimalPresentationMapInv,MinimalPrimes,minimalPrimes,minimalReduction,Minimize,minimizeFilename,MinimumVersion,minors,minPosition,minPres,minprimes,Minus,minus,Miura,MixedMultiplicity,mkdir,mod,Module,module,ModuleDeformations,modulo,MonodromySolver,Monoid,monoid,MonoidElement,Monomial,MonomialAlgebras,monomialCurveIdeal,MonomialIdeal,monomialIdeal,MonomialIntegerPrograms,MonomialOrbits,MonomialOrder,Monomials,monomials,MonomialSize,monomialSubideal,moveFile,multidegree,multidoc,multigraded,MultigradedBettiTally,MultiGradedRationalMap,multiplicity,MultiplicitySequence,MultiplierIdeals,MultiplierIdealsDim2,MultiprojectiveVarieties,mutable,MutableHashTable,mutableIdentity,MutableList,MutableMatrix,mutableMatrix,NAGtypes,Name,nanosleep,Nauty,NautyGraphs,NCAlgebra,NCLex,needs,needsPackage,Net,net,NetFile,netList,new,newClass,newCoordinateSystem,NewFromMethod,newline,NewMethod,newNetFile,NewOfFromMethod,NewOfMethod,newPackage,newRing,nextkey,nextPrime,nil,NNParser,NoetherianOperators,NoetherNormalization,NonAssociativeProduct,NonminimalComplexes,nonspaceAnalyzer,NoPrint,norm,normalCone,Normaliz,NormalToricVarieties,not,Nothing,notify,notImplemented,NTL,null,nullaryMethods,nullhomotopy,nullParser,nullSpace,Number,number,NumberedVerticalList,numcols,numColumns,numerator,numeric,NumericalAlgebraicGeometry,NumericalCertification,NumericalImplicitization,NumericalLinearAlgebra,NumericalSchubertCalculus,numericInterval,NumericSolutions,numgens,numRows,numrows,odd,oeis,of,ofClass,OL,OldPolyhedra,OldToricVectorBundles,on,OneExpression,OnlineLookup,OO,oo,ooo,oooo,openDatabase,openDatabaseOut,openFiles,openIn,openInOut,openListener,OpenMath,openOut,openOutAppend,operatorAttributes,Option,OptionalComponentsPresent,optionalSignParser,Options,options,OptionTable,optP,or,Order,order,OrderedMonoid,orP,OutputDictionary,Outputs,override,pack,Package,package,PackageCitations,PackageDictionary,PackageExports,PackageImports,PackageTemplate,packageTemplate,pad,pager,PairLimit,pairs,PairsRemaining,PARA,Parametrization,parent,Parenthesize,Parser,Parsing,part,Partition,partition,partitions,parts,path,pdim,peek,PencilsOfQuadrics,Permanents,permanents,permutations,pfaffians,PHCpack,PhylogeneticTrees,pi,PieriMaps,pivots,PlaneCurveSingularities,plus,poincare,poincareN,Points,polarize,poly,Polyhedra,Polymake,PolynomialRing,Posets,Position,position,positions,PositivityToricBundles,POSIX,Postfix,Power,power,powermod,PRE,Precision,precision,Prefix,prefixDirectory,prefixPath,preimage,prepend,presentation,pretty,primaryComponent,PrimaryDecomposition,primaryDecomposition,PrimaryTag,PrimitiveElement,Print,print,printerr,printingAccuracy,printingLeadLimit,printingPrecision,printingSeparator,printingTimeLimit,printingTrailLimit,printString,printWidth,processID,Product,product,ProductOrder,profile,profileSummary,Program,programPaths,ProgramRun,Proj,Projective,ProjectiveHilbertPolynomial,projectiveHilbertPolynomial,ProjectiveVariety,promote,protect,Prune,prune,PruneComplex,pruningMap,Pseudocode,pseudocode,pseudoRemainder,Pullback,PushForward,pushForward,Python,QQ,QQParser,QRDecomposition,QthPower,Quasidegrees,QuaternaryQuartics,QuillenSuslin,quit,Quotient,quotient,quotientRemainder,QuotientRing,Radical,radical,RadicalCodim1,radicalContainment,RaiseError,random,RandomCanonicalCurves,RandomComplexes,RandomCurves,RandomCurvesOverVerySmallFiniteFields,RandomGenus14Curves,RandomIdeals,randomKRationalPoint,RandomMonomialIdeals,randomMutableMatrix,RandomObjects,RandomPlaneCurves,RandomPoints,RandomSpaceCurves,Range,rank,RationalMaps,RationalPoints,RationalPoints2,ReactionNetworks,read,readDirectory,readlink,readPackage,RealField,RealFP,realPart,realpath,RealQP,RealQP1,RealRoots,RealRR,RealXD,recursionDepth,recursionLimit,Reduce,reducedRowEchelonForm,reduceHilbert,reductionNumber,ReesAlgebra,reesAlgebra,reesAlgebraIdeal,reesIdeal,References,ReflexivePolytopesDB,regex,regexQuote,registerFinalizer,regSeqInIdeal,Regularity,regularity,relations,RelativeCanonicalResolution,relativizeFilename,Reload,remainder,RemakeAllDocumentation,remove,removeDirectory,removeFile,removeLowestDimension,reorganize,replace,RerunExamples,res,reshape,ResidualIntersections,ResLengthThree,Resolution,resolution,ResolutionsOfStanleyReisnerRings,restart,Result,resultant,Resultants,return,returnCode,Reverse,reverse,RevLex,Right,right,Ring,ring,RingElement,RingFamily,ringFromFractions,RingMap,rootPath,roots,rootURI,rotate,round,rowAdd,RowExpression,rowMult,rowPermute,rowRankProfile,rowSwap,RR,RRi,rsort,run,RunDirectory,RunExamples,RunExternalM2,runHooks,runLengthEncode,runProgram,same,saturate,Saturation,scan,scanKeys,scanLines,scanPairs,scanValues,schedule,schreyerOrder,Schubert,Schubert2,SchurComplexes,SchurFunctors,SchurRings,SCRIPT,scriptCommandLine,ScriptedFunctor,SCSCP,searchPath,sec,sech,SectionRing,SeeAlso,seeParsing,SegreClasses,select,selectInSubring,selectVariables,SelfInitializingType,SemidefiniteProgramming,Seminormalization,separate,SeparateExec,separateRegexp,Sequence,sequence,Serialization,serialNumber,Set,set,setEcho,setGroupID,setIOExclusive,setIOSynchronized,setIOUnSynchronized,setRandomSeed,setup,setupEmacs,sheaf,SheafExpression,sheafExt,sheafHom,SheafOfRings,shield,ShimoyamaYokoyama,short,show,showClassStructure,showHtml,showStructure,showTex,showUserStructure,SimpleDoc,simpleDocFrob,SimplicialComplexes,SimplicialDecomposability,SimplicialPosets,SimplifyFractions,sin,singularLocus,sinh,size,size2,SizeLimit,SkewCommutative,SlackIdeals,sleep,SLnEquivariantMatrices,SLPexpressions,SMALL,smithNormalForm,solve,someTerms,Sort,sort,sortColumns,SortStrategy,source,SourceCode,SourceRing,SPACE,SpaceCurves,SPAN,span,SparseMonomialVectorExpression,SparseResultants,SparseVectorExpression,Spec,SpechtModule,SpecialFanoFourfolds,specialFiber,specialFiberIdeal,SpectralSequences,splice,splitWWW,sqrt,SRdeformations,stack,stacksProject,Standard,standardForm,standardPairs,StartWithOneMinor,stashValue,StatePolytope,StatGraphs,status,stderr,stdio,step,StopBeforeComputation,stopIfError,StopWithMinimalGenerators,Strategy,String,STRONG,StronglyStableIdeals,STYLE,Style,style,SUB,sub,SubalgebraBases,sublists,submatrix,submatrixByDegrees,Subnodes,subquotient,SubringLimit,Subscript,subscript,SUBSECTION,subsets,substitute,substring,subtable,Sugarless,Sum,sum,SumOfTwists,SumsOfSquares,SUP,super,SuperLinearAlgebra,Superscript,superscript,support,SVD,SVDComplexes,switch,SwitchingFields,sylvesterMatrix,Symbol,symbol,SymbolBody,symbolBody,SymbolicPowers,symlinkDirectory,symlinkFile,symmetricAlgebra,symmetricAlgebraIdeal,symmetricKernel,SymmetricPolynomials,symmetricPower,synonym,SYNOPSIS,syz,Syzygies,SyzygyLimit,SyzygyMatrix,SyzygyRows,syzygyScheme,TABLE,Table,table,take,Tally,tally,tan,TangentCone,tangentCone,tangentSheaf,tanh,target,Task,taskResult,TateOnProducts,TD,temporaryFileName,tensor,tensorAssociativity,TensorComplexes,terminalParser,terms,TEST,Test,testExample,testHunekeQuestion,TestIdeals,TestInput,tests,TEX,tex,TeXmacs,texMath,Text,TH,then,Thing,ThinSincereQuivers,ThreadedGB,threadVariable,Threshold,throw,Time,time,times,timing,TITLE,TO,to,TO2,toAbsolutePath,toCC,toDividedPowers,toDual,toExternalString,toField,TOH,toList,toLower,top,top,topCoefficients,Topcom,topComponents,topLevelMode,Tor,TorAlgebra,Toric,ToricInvariants,ToricTopology,ToricVectorBundles,toRR,toRRi,toSequence,toString,TotalPairs,toUpper,TR,trace,transpose,TriangularSets,Tries,Trim,trim,Triplets,Tropical,true,Truncate,truncate,truncateOutput,Truncations,try,TSpreadIdeals,TT,tutorial,Type,TypicalValue,typicalValues,UL,ultimate,unbag,uncurry,Undo,undocumented,uniform,uninstallAllPackages,uninstallPackage,Unique,unique,Units,Unmixed,unsequence,unstack,Up,UpdateOnly,UpperTriangular,URL,urlEncode,Usage,use,UseCachedExampleOutput,UseHilbertFunction,UserMode,userSymbols,UseSyzygies,utf8,utf8check,validate,value,values,Variable,VariableBaseName,Variables,Variety,variety,vars,Vasconcelos,Vector,vector,VectorExpression,VectorFields,VectorGraphics,Verbose,Verbosity,Verify,VersalDeformations,versalEmbedding,Version,version,VerticalList,VerticalSpace,viewHelp,VirtualResolutions,VirtualTally,VisibleList,Visualize,wait,WebApp,wedgeProduct,weightRange,Weights,WeylAlgebra,WeylGroups,when,whichGm,while,width,wikipedia,Wrap,wrap,WrapperType,XML,xor,youngest,zero,ZeroExpression,zeta,ZZ,ZZParser}
}
\lstalias{Macaulay2output}{Macaulay2}